\documentclass[11pt,english]{article}

\usepackage{graphicx,psfrag,epsfig,amsfonts,amssymb,amsmath,babel}

 \usepackage[colorlinks=true]{hyperref}

\hypersetup{urlcolor=blue, citecolor=red}

\newtheorem{theorem}{Theorem}[section]

\newtheorem{definition}[theorem]{Definition}

\begin{document}

\title{Coalitions of pulse-interacting dynamical units}

\author{Eleonora Catsigeras\thanks{Instituto de Matem\'{a}tica y Estad\'{\i}stica Rafael Laguardia (IMERL),
 Fac. Ingenier\'{\i}a,  Universidad de la Rep\'{u}blica,  Uruguay.
 E-mail: eleonora@fing.edu.uy
  EC was partially supported by CSIC of Universidad de la Rep\'{u}blica and ANII of Uruguay.} }

\date{\today}

\maketitle

\begin{abstract}

 We prove that large global systems of interacting (non necessarily similar) dynamical units that are coupled by cooperative impulses, recurrently exhibit the so called \em grand coalition \em, for which all the units arrive to their respective goals simultaneously. We bound from above the waiting time until the first grand coalition appears. Finally, we prove that if besides the units are mutually similar, then the grand coalition is the unique  subset of goal-synchronized units that is recurrently shown by the global dynamics.
\end{abstract}

{\noindent \footnotesize {{\em MSC }2010:  {\  Primary: 37NXX, 92B20; Secondary: 34D06, 05C82, 94A17, 92B25} }}

\section{Introduction}

We study the global dynamics of a network $N$ composed by a large number $m  $ of dynamical units that mutually interact by cooperative (i.e. positive) instantaneous pulses.

One of the most cited    examples of the type of phenomena that we are contributing to explain mathematically along this work, is the large scale synchronization of the flashes of   the fireflies \lq\lq Pteroptyx malaccae\rq\rq: a large number of individuals flash periodically all together after a waiting time, when they meet together on trees, with  neither an external   clock nor privileged individuals mastering the global synchronization \cite{ErmentroutFireflies}.

We are motivated on the study of the dynamics of such global systems   to obtain   abstract and very general mathematical results, that are independent of the concrete formulae governing the dynamics, and require very few hypothesis. They are applicable  in particular  to   models  used in Neuroscience for which  more or less  concrete formulae and hypothesis governing the individual dynamics of the neurons are assumed (see  for instance \cite{Dale2,ErmentroutTerman2010,Izhikevich2007,MassBishop2001,stamov2007}).

The mathematical study of the global dynamics of abstract and general networks composed by mutually interacting units has a large diversity of concrete applications to other sciences and technology. As said above, they are widely used in Neuroscience. They have also applications to
  Engineering, for instance in the design and construction of some systems used in communications \cite{YangChua,GlobalSyncSecureCommun2011}; also to Physics, for instance in the study of systems of light controlled oscillators \cite{LCOrecomMartiCabeza2_2003,CabezaMartiRubidoKahan2011}, and in the research of the evolution of physical lattices of coupled dynamical units of different nature \cite{ChazottesFernandez2005,WangSlotine}. They have other important applications to Biology, for instance in the research of mathematical  models of genetic regulatory networks \cite{arnaud}; to Ecology, in the study of the equilibria of some ecological systems evolving on time \cite{ecologiaMutualInterference3_2010,ecology2009}; to Economy and other social sciences, in the research of coupled networks of different agents, individuals or coalitions of individuals, for instance by means of evolutive Game Theory \cite{miltaich,elvio2009}.

While not interacting with the other units of the network, each unit $i \in \{1, 2, \ldots, m\}$, which we also call \lq\lq cell\rq\rq, evolves governed by two rules  that determine the \lq\lq free dynamics of $i$\rq\rq: the \em relaxation rule \em and the \em update rule, \em which we will precisely define in Subsection \ref{subsectionModel}. While the units are not interacting, the dynamics of the network is the product dynamics of its $m$ units, which evolve independently one from the other. But at certain instants, at least one unit $i$ changes the dynamical rules that govern    the other units $j \neq i$. The instants  when each unit $i$ acts on the others  are exclusively determined by the state $x_i$ of $i$. The pulsed coupling  hypothesis assumes that any action  from $i$  to $j \neq i$  is a   discontinuity jump  in the instantaneous state of the cell  $j$ according to the \em interactions rules  \em which we will precisely define in Subsection \ref{subsectionModel}. 

The   free dynamics rules and   the instantaneous interactions rules, as well as the mathematical results that we obtain from them, generalize to a wide context the particular cases that were studied  for instance in \cite{MirolloStrogatz, Bottani1996,Cessac,jimenez2013,YoPierre}.

 The results that we prove along the paper deal with the spontaneous formation of \em coalitions \em (subsets) of dynamical units during the dynamical evolution of the network, provided that the interactions among the units are all   \lq\lq cooperative\rq\rq \ (i.e. positively signed). Roughly speaking, each coalition is a subset of units that synchronize certain milestones of their respective individual dynamics, which we call goals, and do that spontaneously without any external clock or master unit, infinitely many times in the future. In particular the formation of the so called \em grand coalition \em (i.e. the simultaneous arrival to a certain goal of all the units of the network) is spontaneously and recurrently exhibited from any initial state (Theorem \ref{theorem1}). The synchronization in the grand coalition   was initially proved  in 1992 by Mirollo and Strogatz \cite{MirolloStrogatz}, under  restrictive hypothesis  requiring that the units were identical, the interactions were also   identical, and that the free dynamics of the units were one-dimensional oscillators whose evolution were linear on time.   Later, in 1996,  Bottani \cite{Bottani1996} proved the synchronization of the grand coalition requiring that the units were similar (non necessarily identical), but still one dimensional oscillators although their evolution were not necessarily linear on time. In Theorem \ref{theorem1}  we will generalize the result to any  network of non necessarily similar units with cooperative  interactions that depend on the pair of interacting cells, with general free dynamics of each     unit  $i$, on  any  finite  dimension  (depending on $i$), and such that the cells do  not  necessarily behave as oscillators. The price to pay for such a general result is that the network has to be large enough,  and, unless the units were mutually similar (Theorem \ref{theorem2}),  the grand coalition is not necessarily the unique coalition that is exhibited recurrently in the future.

 Due to the fact that the units may be very different and that the grand coalition is not necessarily the unique coalition that is exhibited in the future, the word \lq\lq synchronization\rq\rq \ in Theorem \ref{theorem1}, if applied, it is not  in its classical meaning (\cite{Pikovsky}). In fact, the orbits of each of the units that recurrently exhibit the grand coalition, are not   synchronized in the strict sense since they do not show the same state for all the instants $t \geq 0$. The states of two or more units may sensibly differ one from the others, at some instants between two consecutive formations of the grand coalition.
 
 On the one hand, the synchronization  in the strict or wide sense,  for  models of pulsed coupled dynamical units, were  up to now  proved   for particular examples in which the free dynamics of each cell  is   governed by a differential equation or a discrete time mapping with \em a concrete formulae. \em For instance, the free dynamics is governed by   affine mapping in \cite{Cessac}, by   linear differential equations in \cite{LCOrecomMartiCabeza2_2003,CabezaMartiRubidoKahan2011}, and by piecewise contracting maps in \cite{WangSlotine} \cite{jimenez2013},\cite{YoPierre} or using known results about piecewise contractions  in \cite{bremont}. In this sense, the novelty of the results here is that their proofs work  independently of the concrete formulation of the free dynamics of the cells. They have almost no hypothesis about   the second term   of the differential equation  governing the free dynamics of each of the cells.
 
 On the other hand, the results along this paper  hold  independently of the dimension of the space $X_i$ where the state of each unit evolves, and   they do  not require the free dynamics of each unit to make it an oscillator. This freedom allows the results to be applied for instance to multidimensional chaotic free dynamics of the cells that recurrently shear certain milestones in the global collective dynamics (\cite{LSYoung2008,coombes2013JournMathBio}).
 
 The paper is organized as follows: in Section \ref{section2} we state the mathematical definitions and theorems to be proved. In Section \ref{sectionProofs} we write the proofs.

\section{Definitions and statements of the results}
\label{section2}

\subsection{Definitions and hypothesis} \label{subsectionModel}

\noindent{\bf The relaxation rule of the free dynamics of $i$:}

The relaxation rule of the free dynamics of the cell $i$ determines the evolution on time $t \geq 0$ of the state $x_i$ on a compact finite-dimensional manifold $X_i$ (whose dimension may depend of $i$). It is defined as the solution of any  differential equation:
\begin{equation} \label{eqnFreeDynamics} \frac{d x_i}{dt} = f_i(x_i), \ \       \ x_i  \in X_i\end{equation}
 satisfying just one condition as follows:
  
  There exists a Lyapunov real function $S_i: X_i \mapsto \mathbb{R}$, which we call \em the satisfaction level \em of $i$, such that:
\begin{equation} \label{eqnSatisfaction} \frac{d S_i(x_i(t))}{dt} = \bigtriangledown S_i (x_i(t)) \cdot f_i (x_i(t)) >v_i >0 \ \ \forall \ t \mbox{ such that } S_i(x_i (t)) < \theta_i,\end{equation} where $\theta_i$ is a positive constant (for each unit $i$) which we call the \em goal  \em of $i$. (In formula (\ref{eqnSatisfaction}) $\bigtriangledown S_i\cdot f_i$ denotes the inner product in the tangent bundle of the manifold $X_i$).

In other words, the free dynamics of $i$ holds at all the instants for which $i$ is uncoupled to the network and its state is unchanged by interferences that may come from outside  $i$. It is described by a finite dimensional variable $x_i$ evolving on time $t$ in such a way that the satisfaction level $S_i(x_i)$, while it does not reach the goal value $\theta_i$, is strictly increasing with   $t$  and its (positive) velocity is bounded away from zero.

\vspace{.3cm}

 \noindent{\bf The update rule of the free dynamics of $i$:}
 
 The update rule is a discontinuity jump in the state $x_i$ of the cell $i$ that is produced whenever the satisfaction variable $S_i(x_i(t))$ reaches (or is larger than) the goal level $\theta_i$. This discontinuity jump instantaneously resets the satisfaction level $S_i(x_i(t))$  to a \lq\lq reset value\rq\rq, which is strictly smaller than $\theta_i$. With no loss of generality, we assume that the reset value is zero (see Figure \ref{figure1}). Precisely:
 \begin{equation} \label{eqnResetRule} S_i(x_i(t_0^-)) \geq \theta_i \ \ \Rightarrow \ \ S_i(x_i(t_0)) = 0,\end{equation}
 where $S_i(x_i(t_0^-))$ denotes $ \lim_{t \rightarrow t_0^-} S_i(x_i(t))$.

 Note that the alternation between the relaxation and update rules of the free dynamics of $i$ will occur while no interferences come from outside $i$ forcing its satisfaction variable to decrease (see Figure \ref{figure1}). Nevertheless, the free evolution $S_i(x_i(t))$ is not necessarily periodic  if $\mbox{dim}(X_i) \geq 2$. In fact,  the set $S_i^{-1}(\{0\}) \subset X_i$ of states with constant null satisfaction may be for instance a curve: there may exist infinitely many points in $X_i$ for which $S_i= 0$. So, each  state $x_i(t)$ obtained from the reset rule $S_i(x_i(t)) = 0$  from the goal $S_i(x_i(t^-)) = \theta_i$, does not necessarily repeat in the future to make the evolution $S_i(x_i(t))$ periodic with an exact time-period. On the contrary, if the set of all the possible reset states $x_i \in S_i^{-1}(\{0\}) $  were finite (this can occur even if $S_i^{-1}(\{0\})$ is infinite), then the free dynamics of $i$ would make it be   periodic, i.e. an oscillator.
 
 \begin{definition}
 {\bf (Spikes)}  \em 
 Taking the name from Neuroscience, we call \em spike \em of the cell $i$ to the discontinuity jump of its satisfaction state from the goal value $\theta_i$ (which in Neuroscience is called \lq\lq threshold level\rq\rq) to its reset value (which is assumed to be zero). Note that the instants when each cell $i$ spikes, while not interacting with the other units of the network, are defined just by the value of its own satisfaction variable.   There is neither an external clock nor a master unit in the network to force a synchronization of the spikes of the many cells of the network.
 \end{definition}

{\begin{figure}
[h]
\begin{center}
\vspace{-.4cm}
\includegraphics[scale=.4]{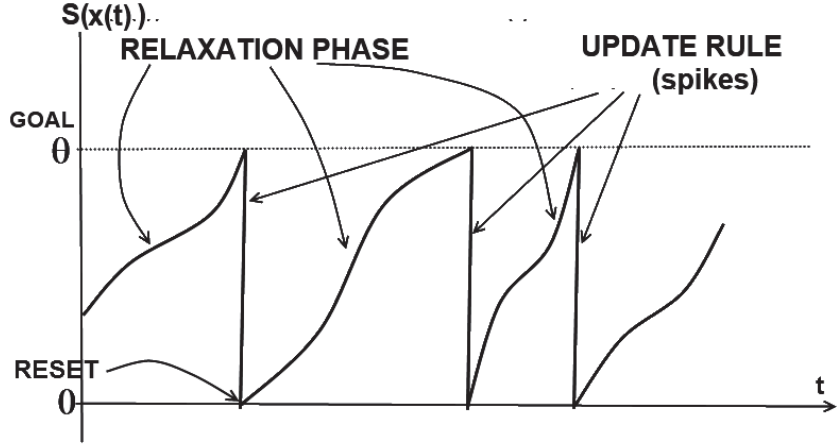}
\vspace{-1.4cm}
\caption{\label{figure1} The evolution on time $t$ of the satisfaction variable $S(x(t))$ of a dynamical unit while not interacting with the other units of the network.}
\vspace{-.5cm}
\end{center}
\end{figure}}

\noindent{\bf The interactions rules among the units} 
 
 Now, let us define the rules that govern the mutual interactions among the units, to compose a global dynamical system which we call network $N$. Consider a system composed by $m \geq 2$ dynamical units with the free dynamics as  described above.

 \begin{definition} {\bf (Spiking instants and inter-spike intervals)}
 \label{definitionSpiking instants} \em We denote by $\{t_n\}_{n \geq 0}$ the sequence of instants $0 \leq t_n < t_{n+1} $ for which at least one cell of the system spikes. We call $t_n$ the \em $n$-th. spiking instant \em of the global system.

 We call $(t_{n+1}, t_n) $ the \em $n$-th. inter-spike interval of the global system. \em
 \end{definition}

 First, by hypothesis, the interactions among the units of the global system are produced only at the spiking instants. In other words, during the inter-spike intervals  the cells evolve independently one from the others.  Hence, the dynamics of the global system  along the inter-spike time intervals is the product dynamics of those of its units.

 Second, at each instant $t_n$ the possible action from a cell $i$ to $j \neq i$ is weighted by a real number $\Delta_{ij}$. The interactions in the network are represented by the edges of a finite graph, whose vertices are the cells $i \in \{1, \ldots, m\}$ and whose edges $(i,j)$ are oriented and weighted by $\Delta_{i,j}$ respectively (see Figure \ref{figure2}). We call $ \Delta_{i,j} $ the interaction  weight. We say that the graph of interactions is \em complete \em if $\Delta_{i,j} \neq 0$ for all $i \neq j$.

 Third and finally, the satisfaction value of any cell $j$, at any spiking instant $t_n$  is defined by the following rule:
 \begin{equation} \label{eqnInteraction}  S_j(x_j(t_n)) = S_j(x_j(t_n^-)) + \sum_{i \in I(t_n), i \neq j} \Delta_{ij} \ \ \ \ \mbox{ if } \ \ \ \ S_j(x_j(t_n^-)) + \sum_{i \in I(t_n), i \neq j} \Delta_{ij} < \theta_j, \end{equation} $$   S_j(x_j(t_n))= 0 \mbox{ otherwise, }$$
 where $I(t_n)$ is the set of neurons that spike at instant $t_n$, and $\Delta_{i,j}$ are the   interactions weights. \begin{definition}\label{definitioncoalition}{\bf (Coalition)}
  
  \em We call the set $I(t_n)$ \em the coalition \em at the spiking instant $t_n$. A coalition $I$ is a singleton if $\#I = 1$. From the definition of the spiking instant, no coalition is empty. \end{definition}

 {\begin{figure}
[h]
\begin{center}
\vspace{-.4cm}
\includegraphics[scale=.4]{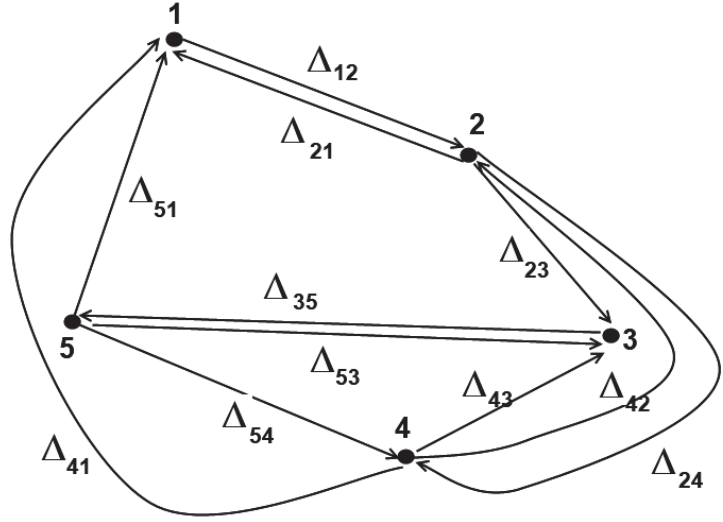}
\vspace{0cm}
\caption{\label{figure2} The graph of interactions of a global system of instantaneously coupled units $1, 2, \ldots, 5$. The oriented and nonzero  weighted edges are denoted by  $\Delta_{ij}$.}
\vspace{-.5cm}
\end{center}
\end{figure}}

  If the interactions weights $\Delta_{i,j} $ are all positive and large enough, the coalition $I(t_n)$ may be the result of  an avalanche process that makes more and more cells spike at the same instant $t_n$ when at least one cell spikes.  In fact, we can always decompose  $I(t_n) $ as the following union of pairwise disjoint (maybe empty) subsets of cells:
  $$I(t_n) = \bigcup_{p \geq 0} I_p(t_n),$$
 where
   $I_0(t_n)$ is the set of cells $i$ such that $x_i(t_n^-)  = \theta_i$, and for all $p \geq 1$, the set $I_{p}(t_n)$ is composed by the cells $j \not \in \cup_{k= 0}^{p-1} I_k(t_n)$ such that $x_j(t_n^-) + \sum_{k= 0}^{k= p-1} \sum_{i \in I_k(t_n)} \Delta_{ij} \geq \theta_{j}. $

   \begin{definition} {\bf (Cooperative and antagonist cells)}
   \label{definitionCooperative} \em

   A cell $i$ is called  \em cooperative \em if $\Delta_{ij} \geq 0$ for all $j \neq i$. It is called \em antagonist \em if $\Delta_{ij} \leq 0$ for all $j \neq i$. It is called \em mixed \em if it is neither cooperative nor antagonist.

   \end{definition}
   In Figure \ref{figure3} we draw the evolution on time of the satisfaction variables of two interacting dynamical units: one of the units is cooperative and the other is antagonist.

   {\begin{figure}
[h]
\begin{center}
\vspace{-.4cm}
\includegraphics[scale=.4]{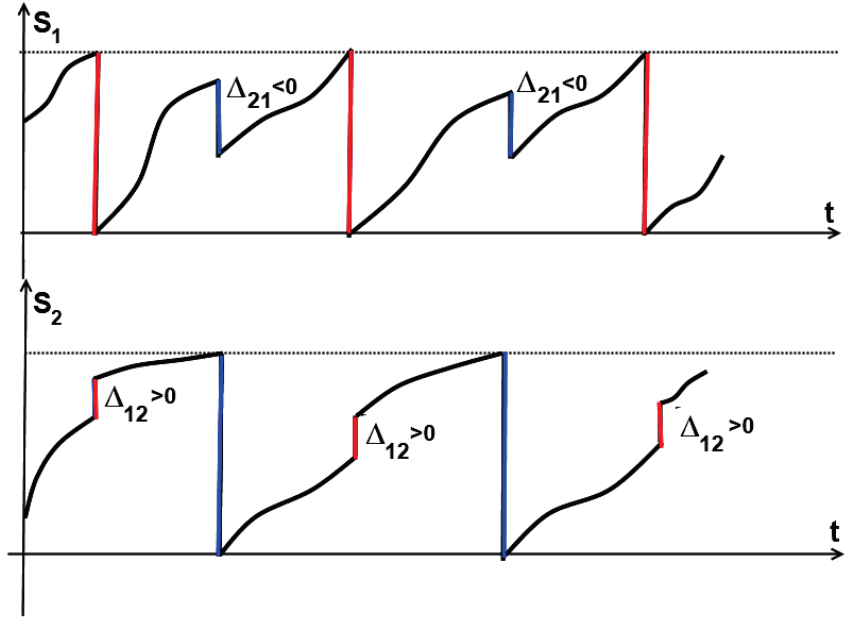}
\vspace{0cm}
\caption{\label{figure3}  Evolution on time $t$  of the satisfaction variable of two interacting units. One cell is cooperative and the other is antagonist.}
\vspace{-.5cm}
\end{center}
\end{figure}}

   From the rule (\ref{eqnInteraction}), when a cooperative cells spikes, it helps the other cells to increase the values of their respective satisfaction variables, so it shortens the time that the others must wait to arrive to their respective goals. On the contrary, an antagonist cell diminishes the values of the satisfaction variables of the other cells, opposing to them and enlarging the time that the others must wait to arrive to their goals.

   Experimentally in Neuroscience,   the nervous system of animals rarely show the existence of mixed cells. This is a reason why one usually assumes the so called Dale's Principle \cite{Dale1, Dale2}:   any cell in the network is either cooperative or antagonist. In \cite{yoDale} abstract mathematical reasons that support Dale's principle were proved: it is a necessary condition for a maximum dynamical richness in the network. Precisely,   the  amount of information that the network can exhibit along its temporal evolution in the future acquires its maximum restricted to a constant number of nonzero interactions, only if Dale's principle holds.

   Along this work we  focuss on the global dynamics of networks that are composed by   cooperative cells and that have a complete graph of interactions.

\vspace{.3cm}

   \noindent{\bf The global state and the vectorial satisfaction variable}
   
   We denote by $${\bf{x}}(t) = (x_1(t), \ldots, x_m(t)) \in \prod_{i= 1}^m X_i$$   the state of the global system at instant $t \geq 0$. We denote  by $${\bf S}({\bf{x}}(t)) = (S_1(x_1(t)), \ldots, S_m(x_m(t))) \in \mathbb{R}^m$$  the vectorial satisfaction variable of the global system at instant $t$. We consider the cube $$Q := \prod_{i= 1}^{m} [0, \theta_i) \subset \mathbb{R}^m.$$ From the hypothesis of the free dynamics of the cells  and of the mutual interactions, if all the cells are cooperative then $${\bf S}({\bf{x}}(t))  \in Q  \ \ \forall \ t \geq 0$$  provided that  \begin{equation} \label{eqnInitialInQ} {\bf{x}}(0) \in {\bf S} ^{-1}(Q).\end{equation} Along this paper we will assume condition (\ref{eqnInitialInQ}). This assumption  is not a restriction for the study of all the orbits of the global autonomous system. In fact, if ${\bf S} ({\bf x}(0)) \not \in Q$, then, applying the inequality (\ref{eqnSatisfaction}) and the reset rule (\ref{eqnResetRule}, and taking into account that the the  interactions are non negative, we deduce that there exists a minimum positive instant $t_0$ such that ${\bf S} ({\bf x}(t_0)) \in Q$. So, translating the origin of the time axis to $t_0$, we have reduced the problem to the case for which the vectorial satisfaction value initially belongs to $Q$.

   \begin{definition} {\bf (Grand coalition)}
   \label{definitionGrandCoalition} \em
   We call $I(t_n)$, defined in \ref{definitioncoalition}, the \em grand coalition \em if all the cells of the system spike at instant $t_n$. Namely, the grand coalition is exhibited at instant $t_n$ if $I(t_n) = \{1, 2, \ldots, m\}$.
   \end{definition}

   \begin{definition} {\bf (Waiting time)}
   \label{definitionWaitingTime} \em If from some initial state of the global system  the grand coalition is exhibited at some spiking instant $t_n \geq 0$, we call the minimum such an instant    \em  the waiting time \em until the grand coalition occurs. Note that  in general,  if existing, the finite waiting time  depends on the initial state.
   \end{definition}
\noindent {\bf Weak interactions:} We will   not need to assume the following condition (\ref{eqnWeakInteractions}) as   an hypothesis. So, it is not an assumption in any part of this paper. Nevertheless, we pose condition (\ref{eqnWeakInteractions})   just because some of the theorems that we will prove along the work    become more interesting for   networks that satisfy it:
\begin{equation}
\label{eqnWeakInteractions}
\max_{i \neq j} |\Delta_{ij}| \ll \min_i \theta_i,\end{equation}
where $\ll$ denotes \lq\lq much smaller than\rq\rq. For instance,  one may be interested in considering $a \ll b$  (where $0< a < b$)  if   $a/b <10^{-3} $. Condition (\ref{eqnWeakInteractions}) says that  the interactions weights are relatively very weak.  

   \begin{definition}
   \label{definitionLargeNetwork} \em {\bf (Large networks)}

   Let $N$ be a   network composed by $m$  cooperative  units, as described above. We say that $N$ is \em large enough \em in relation to  the cooperative interactions  if the following inequality holds:
   \begin{equation}
   \label{eqnLargeInteractions}
    \sqrt m \geq 1 +\frac{\max_i \theta_i}{\min_{i \neq j} \Delta_{ij}}.
   \end{equation}
   Note that,  inequality (\ref{eqnLargeInteractions}) implies that the graph of interactions is complete. In fact   $ \Delta_{ij} \geq 0$ for all $i \neq j$ because the cells are all cooperative, but $$ \Delta_{ij} \neq 0 \ \forall \ i \neq j $$ to make the minimum   in formula (\ref{eqnLargeInteractions}) be nonzero and make this formula hold  for a finite value of $m$.
   \end{definition}

\subsection{Statements of the results}
\label{subsectionStatements}

The   purpose of this paper is to prove the following results:

\begin{theorem}
\label{theorem1} If the network is cooperative and large enough, then from any initial state    the grand coalition is exhibited infinitely many times in the future.

\end{theorem}

\begin{theorem}
\label{theorem1waitingTime}
  If the network is cooperative and large enough, then from any initial state in   ${\bf S}^{-1}(Q)$ the waiting time $t_{n_0}$  before the grand coalition appears for the first time   is upper bounded by:
$$ t_{n_0} \leq \max_{i} \ \frac{\theta_i}{\min_{x_i \in S_i^{-1}[0, \theta_i]} \bigtriangledown S_i(x_i) \cdot f_i(x_i) }. $$

\end{theorem}

\begin{theorem}
\label{theorem2} If the network is cooperative, large enough and if besides all the cells are mutually similar, i.e.  \begin{equation} \label{eqnSimilarCells}
  \frac{\min_i  \big({\theta_i}/{\max_{x_i \in S_i^{-1}[0, \theta_i]}\bigtriangledown S_i(x_i) \cdot f_i(x_i) }\big)}{\max_i  \big({\theta_i}/{\min_{x_i \in S_i^{-1}[0, \theta_i]}\bigtriangledown S_i(x_i) \cdot f_i(x_i) }\big)} \geq 1- \frac{\min_{i \neq j} \Delta_{ij}}{\max_i \theta_i}
 \end{equation} then, from any initial state and after a waiting time  the grand coalition appears at every spiking instant of the system. 
\end{theorem}

Inequality (\ref{eqnSimilarCells}) is satisfied for instance if the cells have mutually identical free dynamics and besides, for each cell $i$,  the maximum and minimum velocities $\bigtriangledown S_i(x_i) \cdot f_i(x_i)$ - according to which the satisfaction variable $S_i$ increases - are not very different. Hypothesis (\ref{eqnSimilarCells})   also holds if the cells   are not identical but   their differences are small enough so the quotient at left in inequality (\ref{eqnSimilarCells}) - which is strictly smaller than 1 - differs from 1 less than $ \displaystyle \frac{\min_{i \neq j} \Delta_{ij}}{\max_i \theta_i}$. If  besides the interactions weights $\Delta_{i,j}$ are much smaller than $\theta_i$ - cf. condition (\ref{eqnWeakInteractions}) -, then the similarity among the cells must be  very notorious  to satisfy the hypothesis of Theorem \ref{theorem2}.

Roughly speaking,  Theorem \ref{theorem2} states that if the cells are similar enough  then, after a waiting time which depends on the initial state of the global system,   the spike of one cell   makes all the other cells also spike at the same instant. In other words, the only recurrent coalition is the grand coalition.

\section{The proofs} \label{sectionProofs}

\subsection{Proof of Theorem \ref{theorem1}}

{\em Proof: } Let $\{t_n\}_{n \geq 0}$ the strictly increasing sequence of spiking instants, as defined in \ref{definitionSpiking instants}. Let
$$ r := 1 + \mbox{int} \Big(\frac{\max_i \theta_i}{\min_{i \neq j} \Delta_{ij}}  \Big ),$$
where int denotes the lower integer part. Since by hypothesis the network is large, from Definition \ref{definitionLargeNetwork} we obtain:
$$r^2 \leq m,$$
where $m$ is the number of units in the system.

It said in Section \ref{section2}, it is not restrictive to assume that the initial state ${\bf x}(0)$ belongs to ${\bf S}^{-1} (Q)$. Thus $S_i(x_i(0)) \in [0, \theta_i)$ for any unit $i$.
We state

\vspace{.3cm}

\noindent{\bf  Assertion (A) } \em During the time interval $[0, t_{r-1}]$ all the units of the system have spiked at least once. \em

\vspace{.3cm}

To prove Assertion (A), let argue by contradiction. Assume that there is a cell, say $j$, such that $x_j(t) < \theta_j$ for all $t \in [0, t_{r-1}]$. By the interactions rule (\ref{eqnInteraction}), and since at least one cell spikes at instant $t_k$ for all $k= 0, \ldots, {r-1}$, we have:
$$S_j(x_j(t_{r-1})) \geq S_j(x_j(0)) + r \, \min_{i \neq j} \Delta_{ij} \geq S_j(x_j(0)) + \theta_j \geq \theta_j,$$
contradicting the initial assumption. So Assertion (A) is proved.

Now, we state

\vspace{.3cm}

\noindent{\bf Assertion (B) } \em  If at some instant $t_n$ at least $r$ cells spike simultaneously, then all the cells spike simultaneously at $t_n$.\em

\vspace{.3cm}

To prove Assertion (B) we have, by hypothesis, $\#I(t_n) \geq r$. Due to the interactions rule  (\ref{eqnInteraction}), for any cell $j \not \in I(t_n) $  we obtain:
$$S_j(x_j(t_n)) \geq S_j(x_j(t_n^-)) + r \, \min_{i \neq j} \Delta_{ij} \geq \theta_j, $$
contradicting the assumption that $j \not \in I(t_n)$. Therefore, all cells are in $I(t_n)$ proving Assertion (B).

\vspace{.3cm}

Consider the $r$ coalitions $I(t_0), I(t_1), \ldots, I(t_{r-1})$. Assertion (A) states that each cell $i$ belongs to at least one of those coalitions. Since the number of different cells is $m \geq r^2$, and the number of coalitions in the above list is $r$, there exists at least one of such coalitions, say  $I(t_k)$ having at least $r$ different cells. In other words, there exists a spiking instant $t_k$ such that at least $r$ cells spike simultaneously at $t_k$. Applying Assertion (B) we deduce that all the cells spike simultaneously at $t_k$. We have proved that  the grand coalition $I(t_k) = \{1, \ldots, m\}$ is spontaneously formed at the instant $t^*_0 := t_k >0$. Since this assertion holds for any initial state, we now translate the origin of the time axis to $t^*_0$, adopting a   new initial state from which the grand coalition will be formed again at some future instant $t^*_ 1 > t^*_0$. By induction on $n$, the grand coalition will be exhibited recurrently in the future at an increasing sequence of instants $t_n^*$, ending the proof of Theorem \ref{theorem1}. \hfill $\Box$

\vspace{.3cm}

\subsection{Proof of Theorem \ref{theorem1waitingTime}}

{\em Proof:}

From the proof of Theorem \ref{theorem1}, the waiting time $t^*_0 $ until the first grand coalition appears is not larger than the instant $t_{r-1}$ such that all the cells have spiked at least once during the time interval $[0, t_{r-1}$. Since all the interactions are positive, $t_{r-1}$ is not larger than the time $T_i$ that the slowest cell, say $i$, would take to arrive to its goal $\theta_i$ if it were not coupled to the network, i.e. under the free dynamics:
 $$t_0^* \leq t_{r-1} \leq T_i.$$ From the relaxation rules (\ref{eqnFreeDynamics}) and (\ref{eqnSatisfaction}) we get
 $$\theta_i = S_i(x_i(T_i^-)) = \int_0^{T_i} \bigtriangledown S_i(x_i(t)) \cdot f_i(x_i(t)) \, dt  \geq \Big(\min_{x_i \in S_i^{-1} ([0, \theta_i])} \bigtriangledown S_i(x_i) \cdot f_i (x_i) \, \Big) \, T_i$$
 Thus 
 $$t_0^* \leq T_i \leq   \frac{\theta_i}{\min_{x_i \in S_i^{-1}[0, \theta_i]} \bigtriangledown S_i(x_i) \cdot f_i(x_i) } \leq \max_{i} \ \frac{\theta_i}{\min_{x_i \in S_i^{-1}[0, \theta_i]} \bigtriangledown S_i(x_i) \cdot f_i(x_i) } ,$$
 ending the proof of Theorem \ref{theorem1waitingTime}. \hfill $\Box$
 
 \vspace{.3cm}
 
 \subsection{Proof of Theorem \ref{theorem2}}
 
 {\em Proof: } 
 
 From Theorem \ref{theorem1}, there exists a first instant $t_0^*$ such that the grand coalition is exhibited. From the update rule (\ref{eqnResetRule}, the state ${\bf x}(t_0^*)$ of the global system is such that ${\bf S}({\bf x}(t_0^*)) = {\bf 0}$. We now translate the origin of the time axis to $t_0^*$. So, the initial state is now ${\bf x}(0)$ such that ${\bf S} ({\bf x}(0)) = {\bf 0}$. 
 
 Hence, to prove Theorem \ref{theorem2} it is enough to show that,   if the hypothesis of inequality (\ref{eqnSimilarCells}) holds, then for any initial state  ${\bf x}(0)$ such that ${\bf S} ({\bf x}(0)) = {\bf 0}$, all the cells spikes simultaneously at the minimum instant $t_1 >0$ such at least one cell, say $i$, spikes.
 
 So, let us compute the values of the satisfaction variables of all the cells at the instant $t_1^-$.
 Due to the relaxation rules (\ref{eqnFreeDynamics}) and (\ref{eqnSatisfaction}) we have
 \begin{equation}
 \label{eqn01}
 S_j(x_j(t_1^-)) = \int _0 ^{t_1} \bigtriangledown S_j(x_j(t)) \cdot f_j(x_j(t)) \, dt   \geq \Big (\min_{x_j \in S_j^{-1}([0, \theta_j]}  \bigtriangledown S_j(x_j) \cdot f_j (x_j)   \Big) \ t_1 \ \  \    \forall \ 1 \leq j \leq m  . \end{equation}
 In particular for the spiking cell $i$ we have
   \begin{equation}
 \label{eqn02}\theta_i = S_i(x_i(t_1^-)) = \int _0 ^{t_1}             \bigtriangledown S_i(x_j(t)) \cdot f_i(x_j(t)) \, dt  \leq \Big(\max_{x_i \in S_i^{-1}([0, \theta_i])} \bigtriangledown S_i(x_i) \cdot f_i (x_i)    \Big) \ t_1.\end{equation}
   Combining inequalities (\ref{eqn01}) and (\ref{eqn02}) we deduce:
   $$S_j(x_j(t_1^-)) \geq  {\theta_i}  \ \frac{\min_{x_j \in S_j^{-1}([0, \theta_j])}  \bigtriangledown S_j(x_j) \cdot f_j (x_j)  }{\max_{x_i \in S_i^{-1}([0, \theta_i])}  \bigtriangledown S_i(x_i) \cdot f_i (x_i)} \geq $$
   $$\geq  \theta_j \ \  \frac{\min_i  \big({\theta_i}/{\max_{x_i \in S_i^{-1}[0, \theta_i]}\bigtriangledown S_i(x_i) \cdot f_i(x_i) }\big)}{\max_j  \big({\theta_j}/{\min_{x_j \in S_j^{-1}[0, \theta_j]}\bigtriangledown S_j(x_j) \cdot f_j(x_j) }\big)} \ \ \ \ \forall \ j \neq i. $$ Using now the hypothesis of inequality (\ref{eqnSimilarCells}, we obtain:
   $$S_j(x_j(t_1^-)) \geq \theta_j \ \Big (1- \frac{\min_{i \neq j} \Delta_{ij}}{\max_i \theta_i}\Big) \geq \theta_j - \min_{i \neq j} \Delta_{ij}  \ \ \ \ \forall \ j \neq i.$$
   Since at least the cell $i$ spikes at instant $t_1$ we have
   $$S_j(x_j(t_1^-)) + \sum_{i \in I(t_1), \ i \neq j} \Delta_{ij} \geq S_j(x_j(t_1^-)) + \min_{i \neq j} \Delta_{ij} \geq \theta_j.$$
   So, applying the interaction rule (\ref{eqnInteraction}) we deduce that the cell $j$ spikes at instant $t_1$. This result holds for all the cells $j \neq i$. Thus, all the cells spike when at least one spikes, ending the proof of Theorem \ref{theorem2}.
\hfill $\Box$

\vspace{.5cm}

\noindent{\bf \large Acknowledgement }

\noindent We thank the scientific and organizing committees of the IV Coloquio Uruguayo de Mate\-m\'{a}\-tica,  for the invitation to give a talk during the event on the subject of this paper.

\end{document}